\documentstyle[amssymb,12pt]{article}

\begin{document}

\newtheorem{theorem}{Theorem}{}
\newtheorem{lemma}[theorem]{Lemma}{}
\newtheorem{corollary}[theorem]{Corollary}{}
\newtheorem{conjecture}[theorem]{Conjecture}{}
\newtheorem{proposition}[theorem]{Proposition}{}
\newtheorem{axiom}{Axiom}{}
\newtheorem{remark}{Remark}{}
\newtheorem{example}{Example}{}
\newtheorem{exercise}{Exercise}{}
\newtheorem{definition}{Definition}{}

\title{A family of combinatorial identities arising form quantum affine 
algebras}
\author{
Alexey Sevostyanov \footnote{e-mail seva@teorfys.uu.se}\\ 
Institute of Theoretical Physics, Uppsala University,\\ and Steklov Mathematical 
Institute, St.Petersburg 
}

\maketitle
\begin{flushright}
UU--ITP 6/98

\end{flushright}

\begin{abstract}
We obtain a family of new combinatorial identities for symmetric formal 
power series.
\end{abstract}

\section*{Introduction}

 This present paper is motivated mainly by the needs of the Drinfeld--Sokolov 
reduction for quantum groups (see \cite{S} for the motivations). However the 
identities which we 
derive are purely combinatorial and may be considered in a more general 
situation. Our main observation is that some solutions of a simple system 
of functional equations satisfy a family of complicated 
combinatorial identities. This phenomenon may lead to new interesting 
effects in the theory of symmetric functions \cite{M}.

\section{A family of combinatorial identities}

 Let $a_{ij} , i,j=1,\ldots , l$ be a generalized Cartan matrix : 
$a_{ii}=2$ , $a_{ij}$ are nonpositive integers for $i \neq j$ and 
$a_{ij}=0$ implies $a_{ji}=0$. Suppose 
also that $a_{ij}$ is symmetrizable , i. e. , there are coprime positive 
integers $d_1,\ldots ,d_l$ such that the matrix $b_{ij}=a_{ij}d_i$ is 
symmetric. 

We shall make use of formal power series ( f.p.s. ) which are infinite in both 
directions. The space of such series is denoted by ${\Bbb C}((z))$. The product 
of 
two f.p.s. $f(z)=\sum_{n=-\infty}^{\infty}f_nz^n , 
g(z)=\sum_{n=-\infty}^{\infty}g_nz^n$ is 
said to exist if the coefficients of the series

\[
\sum_{p=-\infty}^{\infty}z^p\sum_{k+n=p}f_ng_k
\]

are well defined. That is the series

\[
\sum_{k+n=p}f_ng_k
\]

converges for every p.

Similarly , the product of three series $f(z)=\sum_{n=-\infty}^{\infty}f_nz^n ,
g(z)=\sum_{n=-\infty}^{\infty}g_nz^n , h(z)=\sum_{n=-\infty}^{\infty}h_nz^n $ 
exists if the 
series

\[
\sum_{k+n+l=p}f_ng_kh_l
\]

converges for every p and its sum does not depend on the ordering of the terms.
Clearly , in this case the products $g(z)h(z) , g(z)f(z)$ and $f(z)h(z)$ are 
well--defined.

For instance , if two or more series have a common domain of convergence their 
product is well--defined.

Denote by $\frac{1}{1-x}$ the geometric series 

\[
\frac{1}{1-x}=\sum_{n=0}^\infty x^n.
\]

Let $\Sigma_n$ be the symmetric group of order $n$.
For $q \in {\Bbb C}$ we shall consider the following system of equations for 
f.p.s :

\begin{equation}\label{1}
(u-vq^{ b_{ij}})F_{ji}(z^{-1})=(q^{ b_{ij}}u-v)F_{ij}( z) , a_{ij}\neq 0 ,
\end{equation}

\begin{equation}\label{2}
\begin{array}{l}
\sum_{\pi \in \Sigma_{1-a_{ij}}}\sum_{k=0}^{1-a_{ij}}(-1)^k 
\left[ \begin{array}{c} 1-a_{ij} \\ k \end{array} \right]_{q_i}
\prod_{p<q}F_{ii}(\frac {z_{\pi (q)}}{z_{\pi (p)}} )\prod_{r=1}^k F_{ji}(\frac 
{w}{z_{\pi (r)}})\times \\
\prod_{s=k+1}^{1-a_{ij}}F_{ij}(\frac {z_{\pi (s)}}{w})=0 , a_{ij} \neq 0 , i 
\neq j , q_i=q^{d_i} \\
 \\
\mbox{ where } \left[ \begin{array}{c} m \\ n \end{array} \right]_q={[m]_q! 
\over [n]_q![n-m]_q!} , 
[n]_q!=[n]_q\ldots [1]_q , [n]_q={q^n - q^{-n} \over q-q^{-1} }.
\end{array}
\end{equation}

Denote the l.h.s of (\ref{2}) by $P_{ij}(z_1, \ldots, z_{1-a_{ij}},w)$. Let 
$|q|<1$.
Put

\[
s_{ij}(z_1, \ldots, z_{1-a_{ij}})=\prod_{p\neq q}
{1 \over 1-q^{b_{ii}}\frac {z_q}{z_p}},
\]

\[
r_{ij}(z_1, \ldots, z_{1-a_{ij}},w)=\prod_{s=1}^{1-a_{ij}}
{1 \over 1-q^{b_{ij}}\frac {z_s}{w}}.
\]

Our main goal is to prove the following theorem:

\begin{theorem}\label{th1}
Let $F_{kl}(z),k,l=1,\ldots , l$ be a solution of system (\ref{1}). 
Suppose that for some i and j the product

\begin{equation}\label{constr1}
s_{ij}\cdot r_{ij} \cdot \prod_{p\neq q}(1-q^{b_{ii}}\frac 
{z_q}{z_p})\prod_{s=1}^{1-a_{ij}}
(1-q^{b_{ij}}\frac {z_s}{w})P_{ij}
\end{equation}

is well--defined in the sense of f.p.s.. 

Then 

\[
P_{ij}=0 .
\]

in the sense of f.p.s.. Thus the solution $F_{kl}(z)$ satisfies to the identity 
(\ref{2}).

\end{theorem}

We divide the proof of the theorem into several lemmas.

\begin{lemma}\label{l1}
Let $F_{kl}(z),k,l=1,\ldots , l$ be a solution of system (\ref{1}). Then 

\begin{equation}\label{*}
\prod_{p\neq q}(1-q^{b_{ii}}\frac {z_q}{z_p})\prod_{s=1}^{1-a_{ij}}
(1-q^{b_{ij}}\frac {z_s}{w})P_{ij}=0 , i \neq j , a_{ij} \neq 0 .
\end{equation}

\end{lemma}

{\em Proof.} Let $\pi \in \Sigma_{1-a_{ij}}$. Consider the product:

\[
\prod_{s=1}^{1-a_{ij}}(1-q^{b_{ij}}\frac {z_s}{w})\prod_{r=1}^k F_{ji}(\frac 
{w}{z_{\pi (r)}})
\prod_{s=k+1}^{1-a_{ij}}F_{ij}(\frac {z_{\pi (s)}}{w}).
\]

The f.p.s. $\prod_{s=1}^{1-a_{ij}}(1-q^{b_{ij}}\frac {z_s}{w})$ is symmetric 
with 
respect to permutations of the formal variables $z_s$. Therefore

\[
\prod_{s=1}^{1-a_{ij}}(1-q^{b_{ij}}\frac {z_s}{w})\prod_{r=1}^k F_{ji}(\frac 
{w}{z_{\pi (r)}})
\prod_{s=k+1}^{1-a_{ij}}F_{ij}(\frac {z_{\pi (s)}}{w})=
\]

\[
\prod_{s=1}^{1-a_{ij}}(1-q^{b_{ij}}\frac {z_\pi (s)}{w})\prod_{r=1}^k 
F_{ji}(\frac {w}{z_{\pi (r)}})
\prod_{s=k+1}^{1-a_{ij}}F_{ij}(\frac {z_{\pi (s)}}{w}).
\]

Now using equations (\ref{1}) for $F_{ij}$ we obtain:

\[
\prod_{s=1}^{1-a_{ij}}(1-q^{b_{ij}}\frac {z_\pi (s)}{w})\prod_{r=1}^k 
F_{ji}(\frac {w}{z_{\pi (r)}})
\prod_{s=k+1}^{1-a_{ij}}F_{ij}(\frac {z_{\pi (s)}}{w})=
\]

\[
\prod_{s=1}^k(1-q^{b_{ij}}\frac {z_{\pi 
(s)}}{w})\prod_{s=k+1}^{1-a_{ij}}(q^{b_{ij}}-\frac {z_{\pi (s)}}{w})
\prod_{r=1}^{1-a_{ij}}F_{ji}(\frac {w}{z_{\pi (r)}})=
\]

\begin{equation}\label{3}
\prod_{s=1}^k(1-q^{b_{ij}}\frac {z_{\pi 
(s)}}{w})\prod_{s=k+1}^{1-a_{ij}}(q^{b_{ij}}-\frac {z_{\pi (s)}}{w})
\prod_{r=1}^{1-a_{ij}}F_{ji}(\frac {w}{z_{r}}),
\end{equation}

since $\prod_{r=1}^{1-a_{ij}}F_{ji}(\frac {w}{z_{\pi (r)}})$ is also a symmetric 
f.p.s..

Similarly,

\[
\prod_{p\neq q}(1-q^{b_{ii}}\frac {z_q}{z_p})\prod_{p<q}F_{ii}(\frac {z_{\pi 
(q)}}{z_{\pi (p)}} )=
\]

\[
\prod_{p\neq q}(1-q^{b_{ii}}\frac {z_{\pi (q)}}{z_{\pi (p)}})\prod_{p<q}
F_{ii}(\frac {z_{\pi (q)}}{z_{\pi (p)}} )=
\]

\begin{equation}\nonumber
\begin{array}{l}
\prod_{p<q}\left( \prod_{\pi (q)>\pi (p)}F_{ii}(\frac {z_{\pi (q)}}{z_{\pi (p)}} 
) 
\prod_{\pi (q)<\pi (p)}F_{ii}(\frac {z_{\pi (q)}}{z_{\pi (p)}} ) \right) \times 
\\
\prod_{p<q}\left( \prod_{\pi (q)>\pi (p)}(1-q^{b_{ii}}\frac {z_{\pi (q)}}{z_{\pi 
(p)}})
\prod_{\pi (q)<\pi (p)}(1-q^{b_{ii}}\frac {z_{\pi (q)}}{z_{\pi (p)}})\right) 
\times \\
\prod_{p>q}(1-q^{b_{ii}}\frac {z_{\pi (q)}}{z_{\pi (p)}})=
\end{array}
\end{equation}

\begin{equation}\nonumber
\begin{array}{l}
\prod_{p<q}\left( \prod_{\pi (q)>\pi (p)}(q^{b_{ii}}-\frac {z_{\pi (q)}}{z_{\pi 
(p)}})
\prod_{\pi (q)<\pi (p)}(1-q^{b_{ii}}\frac {z_{\pi (q)}}{z_{\pi (p)}})\right) 
\times \\
\prod_{p>q}(1-q^{b_{ii}}\frac {z_{\pi (q)}}{z_{\pi (p)}}) \times \\
\prod_{p<q}\left( \prod_{\pi (q)>\pi (p)}F_{ii}(\frac {z_{\pi (p)}}{z_{\pi (q)}} 
)
\prod_{\pi (q)<\pi (p)}F_{ii}(\frac {z_{\pi (q)}}{z_{\pi (p)}} ) \right)=
\end{array}
\end{equation}

\begin{equation}\label{4}
\begin{array}{l}
\prod_{p>q}\left((1-q^{b_{ii}}\frac {z_q}{z_p})F_{ii}(\frac {z_q}{z_p} ) \right) 
\times \\
\prod_{p>q,\pi (q)>\pi (p)}(1-q^{b_{ii}}\frac {z_{\pi (q)}}{z_{\pi (p)}})
\prod_{p<q,\pi (q)>\pi (p)}(q^{b_{ii}}-\frac {z_{\pi (q)}}{z_{\pi (p)}}).
\end{array}
\end{equation}

Substituting (\ref{3}) and (\ref{4}) into (\ref{*}) we get :

\begin{equation}
\prod_{p\neq q}(1-q^{b_{ii}}\frac {z_q}{z_p})\prod_{s=1}^{1-a_{ij}}
(1-q^{b_{ij}}\frac {z_s}{w})P_{ij}=
\end{equation}

\begin{equation}\nonumber
\begin{array}{l}
\prod_{r=1}^{1-a_{ij}}F_{ji}(\frac {w}{z_{r}})
\prod_{p>q}\left((1-q^{b_{ii}}\frac {z_q}{z_p})F_{ii}(\frac {z_q}{z_p} ) \right) 
\times \\
\sum_{\pi \in \Sigma_{1-a_{ij}}}\sum_{k=0}^{1-a_{ij}}(-1)^k 
\left[ \begin{array}{c} 1-a_{ij} \\ k \end{array} \right]_{q_i} \times \\
\prod_{p>q,\pi (q)>\pi (p)}(1-q^{b_{ii}}\frac {z_{\pi (q)}}{z_{\pi (p)}})
\prod_{p<q,\pi (q)>\pi (p)}(q^{b_{ii}}-\frac {z_{\pi (q)}}{z_{\pi (p)}}) \times 
\\
\prod_{s=1}^k(1-q^{b_{ij}}\frac {z_{\pi (s)}}{w})
\prod_{s=k+1}^{1-a_{ij}}(q^{b_{ij}}-\frac {z_{\pi (s)}}{w}).
\end{array}
\end{equation}

Put $q_i=t$. Lemma \ref{l1} follows from

\begin{lemma}\label{l2}
For any $m\in {\Bbb Z} , m \leq 0$ the following identity holds:

\begin{equation}\label{trel}
\begin{array}{l}
\sum_{\pi \in \Sigma_{1-m}}\sum_{k=0}^{1-m}(-1)^k 
\left[ \begin{array}{c} 1-a_{ij} \\ k \end{array} \right]_{t} \times \\
\prod_{p>q,\pi (q)>\pi (p)}(1-t^2\frac {z_{\pi (q)}}{z_{\pi (p)}})
\prod_{p<q,\pi (q)>\pi (p)}(t^2-\frac {z_{\pi (q)}}{z_{\pi (p)}}) \times \\
\prod_{s=1}^k(1-t^m\frac {z_{\pi (s)}}{w})
\prod_{s=k+1}^{1-a_{ij}}(t^m-\frac {z_{\pi (s)}}{w})=0.
\end{array}
\end{equation}

\end{lemma}

{\em Proof.} (see example (\ref{ex1}) below)

{\em Proof of the theorem.} The conditions of the theorem imply that 

\[
s_{ij}\cdot r_{ij} \cdot \prod_{p\neq q}(1-q^{b_{ii}}\frac 
{z_q}{z_p})\prod_{s=1}^{1-a_{ij}}
(1-q^{b_{ij}}\frac {z_s}{w})P_{ij}=P_{ij}.
\]

But by lemma \ref{l1} 

\[
(1-q^{b_{ii}}\frac {z_q}{z_p})\prod_{s=1}^{1-a_{ij}}
(1-q^{b_{ij}}\frac {z_s}{w})P_{ij}=0.
\]

Therefore

\[
P_{ij}=0.
\]

This concludes the proof.

One can formulate several versions of theorem \ref{th1}. For example , put 

\[
r'_{ij}(z_1, \ldots, z_{1-a_{ij}},w)=\prod_{s=1}^{1-a_{ij}}
{\frac {w}{z_s} \over 1-q^{b_{ij}}\frac {w}{z_s}}.
\]

\begin{theorem}\label{th2}
Let $F_{kl}(z),k,l=1,\ldots , l$ be a solution of system (\ref{1}). 
Suppose that for some i and j the product

\begin{equation}\label{constr2}
s_{ij}\cdot r'_{ij} \cdot \prod_{p\neq q}(1-q^{b_{ii}}\frac 
{z_q}{z_p})\prod_{s=1}^{1-a_{ij}}
(\frac {z_s}{w}-q^{b_{ij}})P_{ij}
\end{equation}

is well--defined in the sense of f.p.s.. 

Then 

\[
P_{ij}=0 .
\]

in the sense of f.p.s.. Thus the solution $F_{kl}(z)$ satisfies the identity 
(\ref{2}).

\end{theorem}

Similar statements exist for $|q|>1$.

\section{Examples of combinatorial identities}

\begin{example}{\bf Polynomial identities}\label{ex1}
\end{example}
First consider the single equation

\begin{equation}\label{c}
(z-c)F(z^{-1},c)=(cz-1)F(z,c)
\end{equation}

\begin{remark}\label{r1}

If $F(z,c)$ is a solution of the equation then 
$\Phi (z,c)=F(z^{-1},c^{-1})$ is also a solution. Thus the transformation 
$(z,c) \mapsto (z^{-1},c^{-1})$ is a symmetry of the equation.

\end{remark}

\begin{lemma}\label{l3}
The elements

\begin{equation}\label{basis}
P_n(z,c)=z^n +cz^{-n} +\sum_{p=1-n}^{-1}c^{p+n-1}(c^2 -1) z^p - c^{n-1} (c+1) , 
n \geq 1
\end{equation}

form a basis in the space of Laurent polynomial solutions of equation (\ref{c}).

\end{lemma}

By remark \ref{r1} the set

\[
Q_n(z,c)=P_n(z^{-1},c^{-1}) ,n \geq 1
\]

is another basis in the space of Laurent polynomial solutions.

{\em Proof.} Let $P(z)=\sum_{-N}^M C_nz^n , N,M >0$ be a solution of (\ref{c}). 
Using 
the equation one can express the coefficients $C_n , n\leq 0$ via $C_n , n>0$. 
Simple 
calculation shows that actually $M=N$ and

\[
P(z)=\sum_{n=1}^NC_nP_n(z,c).
\]

The proof follows.

Put $F_{ij}=P_n(z,q^{b_{ij}})$. This gives a solution of system (\ref{1}). The 
conditions of 
theorem \ref{th1} are satisfied for every $a_{ij}$ since $P_n(z,q^{b_{ij}})$ are 
Laurent 
polynomials. Thus 
we obtain a family of combinatorial identities (\ref{2}) for the elements 
$P_n(z,q^{b_{ij}})$. 

Similarly one derives a family of identities for the elements 
$Q_n(z,q^{b_{ij}})$.

Consider the solution $P_1(z,q^{b_{ij}})$ in more detail. By (\ref{basis}) we 
have:

\[
P_1(z,q^{b_{ij}})=z+q^{b_{ij}}z^{-1} -(q^{b_{ij}}+1)=-(1-z)(1-q^{b_{ij}}z^{-1}).
\]

The identities (\ref{2}) for $P_1(z,q^{b_{ij}})$ amount to the relations:

\begin{equation}\label{5}
\begin{array}{l}
\sum_{\pi \in \Sigma_{1-a_{ij}}}\sum_{k=0}^{1-a_{ij}}(-1)^k 
\left[ \begin{array}{c} 1-a_{ij} \\ k \end{array} \right]_{q_i}
\prod_{p<q}(1-\frac {z_{\pi (q)}}{z_{\pi (p)}} )(1-q^{b_{ii}}\frac {z_{\pi 
(p)}}{z_{\pi (q)}}) \times \\
\prod_{r=1}^k (1-\frac {w}{z_{\pi (r)}})(1-q^{b_{ij}}\frac {z_{\pi 
(r)}}{w})\times \\
\prod_{s=k+1}^{1-a_{ij}}(1-\frac {z_{\pi (s)}}{w})(1-q^{b_{ij}}\frac {w}{z_{\pi 
(s)}})=0 , a_{ij} \neq 0 , i \neq j .
\end{array}
\end{equation}

Pulling out the symmetric factor $\prod_{r=1}^{1-a_{ij}} (1-\frac {w}{z_{\pi 
(r)}})$ 
and the antisymmetric one $\prod_{p<q}(\frac {1}{z_{\pi (q)}}-\frac {1}{z_{\pi 
(p)}} )$ 
we find that the identities (\ref{5}) are equivalent to the following relations:

\begin{equation}\label{6}
\begin{array}{l}
\sum_{\pi \in \Sigma_{1-m}}(-1)^{l(\pi )}\sum_{k=0}^{1-m}(-1)^k 
\left[ \begin{array}{c} 1-m \\ k \end{array} \right]_{t}
\prod_{p<q}({z_{\pi (q)}}-t^2{z_{\pi (p)}} ) \times \\
\prod_{r=1}^k (1-t^m\frac {z_{\pi (r)}}{w})
\prod_{s=k+1}^{1-m}(\frac {z_{\pi (s)}}{w}-t^m)=0 , m\in {\Bbb Z} , m \leq 0 ,
\end{array}
\end{equation}

which coincide with the identities obtained by Jing in \cite{J}. Since the 
solution 
$P_1(z,q^{b_{ij}})$ is formed by Laurent polynomials the identity (\ref{6}) is 
equivalent to the relation (\ref{trel}). In particular, this proves lemma 
\ref{l2}.

\begin{example}{\bf Identities for Taylor series}
\end{example}

Now let $F_{ij} \in {\Bbb C}[[z]]$.

\begin{lemma}
The system of equations (\ref{1}) 
has a unique nontrivial solution 
in ${\Bbb C}[[z]]$ with the asymptotics $F_{ij}(0)=q_j^{n_{ij}}$. The solution 
has the form:

\begin{equation}\label{F}
F_{ij}( z)={q_j^{n_{ij}} -zq_i^{n_{ji}} \over 1-zq^{b_{ij}}}, a_{ij}\neq 0 ,
\end{equation}

\end{lemma}

{\em Proof.}  Put

\[
F_{ij}(z)=\sum_{n=0}^\infty c_{ij}^nz^n .
\]

The l.h.s. of (\ref{1}) must belong to ${\Bbb C}[[z]]$. This allows us to 
determine 
$F_{ij}(z)$ up tp a constant:

\[
F_{ij}(z)=c_{ij}^0+c_{ij}^1{z \over 1-zq^{b_{ij}}} , c_{ij}^0=q_j^{n_{ij}} .
\]

Substituting this ansatz into (\ref{1}) we get the following relation for the 
coefficients $c_{ij}^0,c_{ij}^1$:

\[
c_{ij}^1=-c_{ji}^0+q^{b_{ij}}c_{ij}^0 .
\]

This yields (\ref{F}). 

\begin{theorem}
For every $n_{ij}$ such that $d_in_{ji}-d_jn_{ij}=\varepsilon _{ij}b_{ij} , 
\varepsilon _{ij}=-\varepsilon _{ji} , \varepsilon _{ij}=\pm 1$ 
the solution (\ref{F}) satisfies the identities (\ref{2}).
\end{theorem}

{\em Proof} An important property of the solution (\ref{F}) subject to the 
conditions of 
the theorem is that either 
$F_{ji}=q_i^{n_{ji}}$ or $F_{ij}=q_j^{n_{ij}}$. Using this fact one can show 
that 
either the series in (\ref{constr1}) or the series in (\ref{constr2}) have a 
common domain of 
convergence. This allows us to apply theorem \ref{th1} or theorem \ref{th2} 
,respectively, to 
obtain identities (\ref{2}) for the solution (\ref{F}).

\begin{example}{\bf Identities for general formal power series}
\end{example}

\begin{lemma}
The elements

\[
\varphi_n(z,c)=z^{-n}+z^n{c-z \over 1-cz} , n \geq 1
\]

\[
\varphi_0(z,c)={1-z \over 1-cz}
\]

are solutions of equation (\ref{c}). Every solution of equation (\ref{c}) may be 
uniquely
represented as a Taylor series $\sum_{n=0}^\infty C_n \varphi_n(z,c)$.

\end{lemma}

{\em Proof} is quite similar to that of lemma \ref{l3} .

We shall say that the set $\{ \varphi_n(z,c) \}_{n \geq 0}$ is a basis in the 
space of 
f.p.s. solutions of (\ref{c}). As a consequence we obtain that the series

\[
\psi_n(z,c)=z^{n}+z^{-n}{c^{-1}-z^{-1} \over 1-c^{-1}z^{-1}} , n \geq 1
\]

\[
\psi_0(z,c)={1-z^{-1} \over 1-c^{-1}z^{-1}}
\]

form another basis in the same space.

Thus the sets $\{ \varphi_n(z,q^{b_{ii}}) \}_{n \geq 0}$ and 
$\{ \psi_n(z,q^{b_{ii}}) \}_{n \geq 0}$ are two bases in the space of f.p.s. 
solutions of 
equations (\ref{1}) for $i=j$. 

Let $i \neq j$. 

\begin{lemma}
Let $F_{ij}(z)=\sum_{n=-\infty}^\infty c_{ij}^n z^n , i<j$ be arbitrary f.p.s.. 
Put 

$$
\begin{array}{l}
F_{ji}(z)=C\delta (q^{b_{ij}}z) -C_{ij}^0\left( {z \over 1-q^{b_{ij}}z} + 
{z^{-1} \over 
1-q^{-b_{ij}}z^{-1}}\right) + \\
\sum_{n=1}^\infty \left( C_{ij}^{-n}z^n{q^{b_{ij}}-z \over 1-q^{b_{ij}}z} +
C_{ij}^{n}z^{-n}{q^{-b_{ij}}-z^{-1} \over 1-q^{-b_{ij}}z^{-1}}\right) .
\end{array}
$$

Then the set $\{ F_{ij}(z) , F_{ji}(z) , i< j\}$ is a solution of equations 
(\ref{1}) for 
$i \neq j$. Moreover , every solution may be uniquely represented in this form 
for some 
$C , C_{ij}^n , i< j , n \in {\Bbb Z}$.

\end{lemma}

The lemma can be proved similarly to lemma \ref{l3}.

Thus we have described the space of solutions of the system (\ref{1}) 
completely. 
This yields several examples of identities of the type (\ref{2}). 

For instance , put 

$$
\begin{array}{l}
F_{ii}(z)=\varphi_n(z,q^{b_{ii}}) , \\
F_{ij}(z)=\sum_{n=-N}^{-M }C_{ij}^n z^n , i< j , N,M>0 ,N>M ,\\
C=0.
\end{array}
$$

Then

\[
F_{ji}(z)=
\sum_{n=M}^{N} C_{ij}^{-n}z^n{q^{b_{ij}}-z \over 1-q^{b_{ij}}z}.
\]

The conditions of theorem \ref{th2} are satisfied for $i < j$ since in that case 
the series in 
(\ref{constr2}) have a common domain of convergence. Hence we obtain a family of 
identities 
(\ref{2}) for the solution.


\begin{thebibliography}{99}

\bibitem{J} N. Jing , {\em Quantum Kac--Moody algebras and vertex 
representations}, 
q-alg/9802036 .

\bibitem{M} I. G. Macdonald , {\em Symmetric functions and Hall polynomials} , 
2nd 
edition , Claredon Press , Oxford , 1995.

\bibitem{S} A. Sevostyanov , Drinfeld--Sokolov reduction for quantum groups , 
math.QA/9805133


\end{thebibliography}
\end{document}